\documentclass[12pt]{article}
\usepackage{graphicx}
\usepackage{color}
\catcode`\@=11 \@addtoreset{equation}{section}

\catcode`\@=12

\newtheorem{Theorem}{Theorem}[section]
\newtheorem{Proposition}{Proposition}[section]
\newtheorem{Lemma}{Lemma}[section]
\newtheorem{Corollary}{Corollary}[section]
\newtheorem{Remark}{Remark}[section]

\newcommand{\bTheorem}[1]{
\begin{Theorem} \label{T#1} }
\newcommand{\eT}{\end{Theorem}}

\newcommand{\bProposition}[1]{
\begin{Proposition} \label{P#1}}
\newcommand{\eP}{\end{Proposition}}

\newcommand{\bLemma}[1]{
\begin{Lemma} \label{L#1} }
\newcommand{\eL}{\end{Lemma}}

\newcommand{\bCorollary}[1]{
\begin{Corollary} \label{C#1} }
\newcommand{\eC}{\end{Corollary}}

\newcommand{\bRemark}[1]{
\begin{Remark} \label{R#1} }
\newcommand{\eR}{\end{Remark}}

\newcommand{\bFormula}[1]{
\begin{equation} \label{#1}}
\newcommand{\eF}{\end{equation}}

\newcommand{\Ov}[1]{\overline{#1}}

\newcommand{\vr}{\varrho}

\newcommand{\vu}{\vc{u}}
\newcommand{\vc}[1]{{\bf #1}}
\newcommand{\qed}{\bigskip \rightline {Q.E.D.} \bigskip}
\newcommand{\Div}{{\rm div}_x}
\newcommand{\Grad}{\nabla_x}

\newcommand{\tn}[1]{\mbox {\F #1}}
\newcommand{\dx}{{\rm d} {x}}
\newcommand{\dt}{{\rm d} t }

\newcommand{\intO}[1]{\int_{\Omega} #1 \ \dx}

\newcommand{\bProof}{{\bf Proof: }}

\font\F=msbm10 scaled 1000

\definecolor{Cgrey}{rgb}{0.85,0.85,0.85}
\definecolor{Cblue}{rgb}{0.50,0.85,0.85}
\definecolor{Cred}{rgb}{1,0,0}
\definecolor{fancy}{rgb}{0.10,0.85,0.10}

\newcommand\Cbox[2]{%
    \newbox\contentbox%
    \newbox\bkgdbox%
    \setbox\contentbox\hbox to \hsize{%
        \vtop{
            \kern\columnsep
            \hbox to \hsize{%
                \kern\columnsep%
                \advance\hsize by -2\columnsep%
                \setlength{\textwidth}{\hsize}%
                \vbox{
                    \parskip=\baselineskip
                    \parindent=0bp
                    #2
                }%
                \kern\columnsep%
            }%
            \kern\columnsep%
        }%
    }%
    \setbox\bkgdbox\vbox{
        \color{#1}
        \hrule width  \wd\contentbox %
               height \ht\contentbox %
               depth  \dp\contentbox
        \color{black}
    }%
    \wd\bkgdbox=0bp%
    \vbox{\hbox to \hsize{\box\bkgdbox\box\contentbox}}%
    \vskip\baselineskip%
}

\date{}

\begin{document}


\title{Uniqueness of rarefaction waves in multidimensional compressible Euler system}

\author{Eduard Feireisl \thanks{The research of E.F. leading to these results has received funding from the European Research Council under the European Union's Seventh Framework
Programme (FP7/2007-2013)/ ERC Grant Agreement 320078} \and Ond\v rej Kreml\thanks{O.K. acknowledges the support of the GA\v CR (Czech Science Foundation) project P201-13-00522S in the general framework of RVO: 67985840.}}

\maketitle

\bigskip

\centerline{Institute of Mathematics of the Academy of Sciences of the Czech Republic}

\centerline{\v Zitn\' a 25, 115 67 Praha 1, Czech Republic}






\maketitle

\bigskip





\begin{abstract}

We show that 1-D rarefaction wave solutions are unique in the class of bounded entropy solutions to the multidimensional compressible Euler system.
Such a result may be viewed as a counterpart of the recent examples of \emph{non-uniqueness} of the shock wave solutions to the Riemann problem, where infinitely many solutions are constructed by the method of convex integration.

\end{abstract}



\section{Introduction}
\label{i}

We consider the compressible Euler system
\bFormula{i1}
\partial_t \vr + \Div (\vr \vu) = 0,
\eF
\bFormula{i2}
\partial_t (\vr \vu) + \Div (\vr \vu \otimes \vu) + \Grad p(\vr) = 0
\eF
describing the time evolution of the density $\vr = \vr(t,x)$ and the velocity $\vu = \vu(t,x)$ of a compressible barotropic fluid.
The problem is supplemented with the Riemann type initial data
\bFormula{i1a}
\vr(0,\cdot) = \vr_0 = \left\{ \begin{array}{c} \vr_L \ \mbox{for}\ x_1 \leq 0, \\  \vr_R \ \mbox{for}\ x_1 > 0, \end{array} \right.
\eF
\bFormula{i2a}
u^1(0, \cdot) = u^1_0 =  \left\{ \begin{array}{c} u^1_L \ \mbox{for}\ x_1 \leq 0, \\  u^1_R \ \mbox{for}\ x_1 > 0, \end{array} \right.
u^j(0, \cdot) = u^j_0 = 0 \ \mbox{for}\ j > 1.
\eF
For the sake of simplicity, we consider the problem in the $2D$-setting although the result holds true in any space dimension.

In order to identify a class of physically \emph{admissible} solutions,
the system (\ref{i1}), (\ref{i2}) is usually augmented by the {energy inequality}
\bFormula{i3}
\partial_t \left( \frac{1}{2} \vr |\vu|^2 + H(\vr) \right) + \Div \left[ \left( \frac{1}{2} \vr |\vu|^2 + H(\vr) \right) \vu \right] +
\Div (p(\vr) \vu ) \leq 0,
\eF
\[
H(\vr) = \vr \int_1^\vr \frac{p(z)}{z^2} \ {\rm d}z.
\]
As is well known, the \emph{Riemann problem} (\ref{i1}-\ref{i3}) admits a mono-dimensional self-similar solution $\vr = \vr( x_1 /t )$,
$u^1 = u^1(x_1/t)$, $u^j = 0$, $j > 1$. Moreover, any such solution consists of at most three constant states connected by shocks or rarefaction
waves. We focus on the situation when the solution contains only rarefaction waves, meaning it is locally Lipschitz for $t > 0$.

Our main goal is to show that the mono-dimensional self-similar Lipschitz solution is unique in the class of all bounded weak solutions to the
$2D-$problem emanating from the same Riemann data. The question is not academic. Recently, Chiodaroli et al. \cite{ChiDelKre}, \cite{ChiKre} showed that uniqueness fails in the
case of solutions containing shocks. In particular, there exist infinitely many admissible solutions satisfying the energy inequality (\ref{i3}).
These solutions develop oscillations in the second component of the velocity field and apparently do not belong to the class $BV$.

Although well-posedness of the Riemann problem in the class of $BV$-solutions is quite well understood (see Bianchini and Bressan \cite{BiaBre}, Chen and Frid \cite{CheFr2}, Chen et al. \cite{CheFr1}, LeFloch \cite{LeFl}) much less seems to be known concerning \emph{bounded} weak solutions. Leger and Vasseur \cite{LegVas} addressed the problem
in the $1D$-setting extending uniqueness of the shock wave solutions to the class of bounded solutions enjoying certain trace property. They use the method of relative entropies proposed by Dafermos \cite{Daf4}. The question of uniqueness for the Riemann problem in the class of \emph{bounded} weak solutions
remains open.

Here, we use the abstract form of the relative entropy inequality (see \cite{FeNoJi}) adapted to problems with \emph{boundary} conditions. The crucial observation
is that the velocity component $u^1$ corresponding to a rarefaction wave must be monotone, more specifically, nondecreasing. As a consequence, all uncontrollable terms in the relative entropy inequality possess a sign providing stability of the Riemann solution.

The paper is organized as follows. In Section \ref{m}, we recall the basic concepts concerning weak solutions to the compressible system and their normal traces, then we state our main result. Section \ref{r} is devoted to the relative entropy inequality and its necessary modifications to accommodate the boundary conditions and singular ``test'' functions. The proof of the main theorem is finished in Section \ref{p}.

\section{Preliminaries, main result}
\label{m}

For definiteness, we consider the problem (\ref{i1}-\ref{i3}), supplemented with the periodic boundary conditions in the $x_2$ variable. Accordingly,
the relevant spatial domain is
\[
\Omega = (- a, a) \times \mathcal{T}^1,
\]
where $\mathcal{T}^1$ denotes the $1D$ (flat) torus, and
where $a > 0$ is a sufficiently large positive number to accommodate the far field conditions for the Riemann problem on a time interval $(0,T)$.

We consider the class of weak solutions that coincide with the Riemann solution outside the interval $(-a,a)$. More specifically, we prescribe the
normal trace of a solution $[\vr, \vu]$ as follows:

\begin{itemize}
\item {\bf Initial state:}
\[
\vr(0,x_1,x_2) = \vr_0(x_1), \ \vr \vu (0, x_1, x_2 ) = \vr_0 \vu_0 (x_1).
\]

\item {\bf Boundary fluxes:}
\[
\vr u^1 (t,-a,x_2) = \vr_L u^1_L,\ \vr u^1 (t,a,x_2) = \vr_R u^1_R;
\]
\[
\left\{ \begin{array}{c}
\left( \vr u^j u^1 + p(\vr) \right) (t, - a,x_2) = \left( \vr_L u^j_L u^1_L + p(\vr_L) \right) , \\ \\ \left(
\vr u^j u^1 +p(\vr) \right) (t, a, x_2) = \left( \vr_R u^j_R u^1_R + p(\vr_R) \right) \end{array} \right\} , \ j=1,2;
\]
\[
\left\{
\begin{array}{c}
\left( \frac{1}{2} \vr |\vu|^2 + H(\vr) + p(\vr) \right) u^1 (t, - a, x_2) =
\left( \frac{1}{2} \vr_L |\vu_L|^2 + H(\vr_L) + p(\vr_L) \right) u^1_L, \\ \\
\left( \frac{1}{2} \vr |\vu|^2 + H(\vr) + p(\vr) \right) u^1 (t, a, x_2) =
\left( \frac{1}{2} \vr_R |\vu_R|^2 + H(\vr_R) + p(\vr_R) \right) u^1_R;
\end{array}
\right\}
\]
\end{itemize}

\noindent where we have set $u^2_L = u^2_R = 0$.

\subsection{Weak formulation}

In agreement with our choice of initial and boundary data, the weak formulation of the problem (\ref{i1}-\ref{i3}) reads as follows:

\begin{itemize}
\item
{\bf Equation of continuity:}

\bFormula{w1}
\intO{ \left[ \vr(\tau,x) \varphi (\tau, x) - \vr_0 (x) \varphi (0, x) \right] }
\eF
\[
+
\int_0^\tau \int_{\mathcal{T}^1}  \vr_R u^1_R \varphi (t, a, x_2)  \ {\rm d}x_2 \ \dt
- \int_0^\tau \int_{\mathcal{T}^1}  \vr_L u^1_L \varphi (t, -a, x_2)  \ {\rm d}x_2 \ \dt
\]
\[
= \int_0^\tau \intO{ \left[ \vr(t,x) \partial_t \varphi(t,x) + \vr \vu (t,x) \cdot \Grad \varphi (t,x) \right] } \ \dt
\]
for any $0 \leq \tau \leq T$, and any test function $\varphi \in C^1([0,T] \times \Ov{\Omega})$.

\item
{\bf Momentum equation:}
\bFormula{w2}
\intO{ \left[ \vr \vu (\tau, x) \cdot \varphi (\tau, x) - \vr_0 \vu_0 (x) \cdot \varphi (0, x) \right] }
\eF
\[
+
\int_0^\tau \int_{\mathcal{T}^1}    \vr_R u^1_R \vu_R \cdot \varphi ( t, a, x_2 )   \ {\rm d}x_2 \ \dt
- \int_0^\tau \int_{\mathcal{T}^1}   \vr_L u^1_L \vu_L \cdot \varphi (t, -a, x_2)  \ {\rm d}x_2 \ \dt
\]
\[
+
\int_0^\tau \int_{\mathcal{T}^1} p(\vr_R) \varphi^1 (t, a, x_2)\ {\rm d}x_2 \ \dt - \int_0^\tau \int_{\mathcal{T}^1} p(\vr_L) \varphi^1 (t, - a, x_2)\ {\rm d}x_2 \ \dt
\]
\[
= \int_0^\tau \intO{ \left[ \vr \vu (t,x) \cdot \partial_t \varphi (t,x) + \vr [\vu \otimes \vu] (t,x) : \Grad \varphi(t,x) + p(\vr)(t,x) \Div \varphi
(t,x)\right] } \ \dt
\]
for any $0 \leq \tau \leq T$, and any $\varphi \in C^1([0,T] \times \Ov{\Omega};R^2)$.

\item
{\bf Energy inequality:}

\bFormula{w3}
\intO{ \left[ \frac{1}{2} \vr |\vu|^2 + H(\vr) \right](\tau, x) } - \intO{ \left[ \frac{1}{2} \vr_0 |\vu_0|^2 + H(\vr_0) \right] }
\eF
\[
+\int_0^\tau \int_{\mathcal{T}^1}  \left[ \frac{1}{2} \vr_R |\vu_R |^2 + H(\vr_R) + p(\vr_R) \right] u^1_R \ {\rm d}x_2 \ \dt -
\int_0^\tau \int_{\mathcal{T}^1}  \left[ \frac{1}{2} \vr_L |\vu_L |^2 + H(\vr_L) + p(\vr_L) \right] u^1_L \ {\rm d}x_2 \ \dt
\]
\[
\leq 0
\]

\end{itemize}

\bRemark{P1}

Note that bounded weak solutions to systems of conservation laws considered may be viewed as $L^\infty$ functions with divergence-measure, in particular, the normal traces
on the boundary of the space-time cylinder $(0,T) \times \Omega$ are well defined bounded measurable functions, see Chen and Frid \cite{ChenFr}, Chen, Torres and Ziemer \cite{ChToZi}.

\eR

\subsection{Main result}

Assume that
\bFormula{pre}
p \in C^1(0, \infty) \cap [0, \infty), \ p(0) = 0, \ p'(\vr) > 0 \ \mbox{for all}\ \vr > 0, \ p \ \mbox{convex in}\ [0, \infty).
\eF
Our main goal is to prove the following result:

\Cbox{Cgrey}{

\bTheorem{m1}
Let the pressure $p$ satisfy (\ref{pre}). Let  $\tilde \vr = \tilde \vr (x_1/t)$, $\tilde \vu = [\tilde u^1(x_1/t),0]$
be the self-similar solution to the Riemann problem consisting of rarefaction waves (locally Lipschitz for $t > 0$) and such that
\bFormula{condition}
{{\rm ess} \inf}_{(0,t) \times R } \tilde \vr > 0.
\eF
Let $[\vr, \vu]$ be a bounded admissible weak solution satisfying (\ref{w1}-\ref{w3}) and such that
\[
\vr \geq 0 \ \mbox{a.a. in}\ (0,T) \times \Omega.
\]

Then
\[
\vr \equiv \tilde \vr, \ \vu \equiv \tilde \vu \ \mbox{in}\ (0,T) \times \Omega.
\]

\eT
}

\bRemark{m1}

Here, the solution $[\tilde \vr, \tilde \vu]$ is extended as constant (periodic) with respect to the $x_2$ variable.
\eR

\bRemark{m2}
Note that according to \cite[Lemma 2.4]{ChiKre} the self-similar solution to the Riemann problem (\ref{i1}--\ref{i3}) consists only of rarefaction waves and satisfies (\ref{condition}) if and only if the initial Riemann data satisfy
\bFormula{Rdata}
\left|\int_{\rho_L}^{\rho_R} \frac{\sqrt{p'(\tau)}}{\tau} {\rm d}\tau\right| \leq u_R^1 - u_L^1 < \int_0^{\rho_L} \frac{\sqrt{p'(\tau)}}{\tau} {\rm d}\tau + \int_0^{\rho_R} \frac{\sqrt{p'(\tau)}}{\tau} {\rm d}\tau.
\eF
\eR

The rest of the paper is devoted to the proof of Theorem \ref{Tm1}.

\section{Relative entropy inequality}
\label{r}

Following Dafermos \cite{Daf4} (see also Berthelin and Vasseur \cite{BeVa1}, Desjardins \cite{DES2}, Leger and Vasseur \cite{LegVas}, among others) we introduce
the \emph{relative entropy} functional in the form
\[
\mathcal{E} \left( \vr, \vu \Big| r, \vc{U} \right) = \frac{1}{2} \vr |\vu - \vc{U}|^2 + \left(H(\vr) - H'(r)(\vr - r) - H(r)\right).
\]
Following the strategy of \cite{FeNoJi} we derive a relative entropy inequality in the general situation when $[\vr, \vu]$ is a (bounded) weak solution of
the compressible Euler system specified through (\ref{w1}-\ref{w3}), while $[r, \vc{U}]$ is an arbitrary pair of \emph{test functions} that are continuously
differentiable in $[0,T] \times \Ov{\Omega}$, and $r > 0$.

\bProposition{r1}

Let $[\vr, \vu]$ be a bounded admissible solution satisfying (\ref{w1}-\ref{w3}) in $(0,T) \times \Omega$, and let $[r, \vc{U}]$ be a pair of functions
such that
\[
r \in C^1([0,T] \times \Ov{\Omega}), \ \vc{U} \in C^1([0,T] \times \Ov{\Omega};R^2), \ r > 0.
\]

Then the following relative entropy inequality
\bFormula{REI}
\intO{ \mathcal{E} \left( \vr, \vu \Big| r, \vc{U} \right)(\tau,x) } - \intO{ \mathcal{E} \left( \vr_0, \vu_0 \Big| r(0,x), \vc{U}(0,x) \right) }
\eF
\[
+ \int_0^\tau \int_{\mathcal{T}^1} \mathcal{E} \left( \vr_R, \vu_R \Big| r (t,a, x_2), \vc{U}(t,a,x_2) \right) u^1_R {\rm d}x_2 \ \dt
\]
\[
-
\int_0^\tau \int_{\mathcal{T}^1} \mathcal{E} \left( \vr_L, \vu_L \Big| r (t,-a, x_2), \vc{U}(t,-a,x_2) \right) u^1_L {\rm d}x_2 \ \dt
\]
\[
+ \int_0^\tau \int_{\mathcal{T}^1} \Big( p(\vr_R) - p(r)(t, a, x_2) \Big)  \left( u^1_R - U^1(t, a, x_2) \right) {\rm d}x_2 \ \dt
\]
\[
- \int_0^\tau \int_{\mathcal{T}^1} \Big( p(\vr_L) - p(r)(t, - a, x_2) \Big) \left( u^1_L - U^1(1, -a, x_2) \right) {\rm d}x_2 \ \dt
\]
\[
 \leq \int_0^\tau \intO{ \left[ \vr \left( \partial_t \vc{U} + \vu \cdot \Grad \vc{U} \right) \cdot (\vc{U} - \vu) + \Big( p(r) -  p(\vr) \Big) \Div \vc{U} \right](t,x) } \ \dt
\]
\[
+ \int_0^\tau \intO{ \Big[ (r - \vr) \partial_t H'(r) + (r \vc{U} - \vr \vu) \cdot \Grad H'(r) \Big] (t,x) } \ \dt
\]
holds for a.a. $\tau \in (0,T)$.

\eP

\bProof

{\bf Step 1}

Using $\varphi = \vc{U}$ as a test function in (\ref{w2}) we obtain
\bFormula{w4}
\intO{ \left[ \vr \vu (\tau, x) \cdot \vc{U} (\tau, x) - \vr_0 \vu_0 \cdot \vc{U} (0, x) \right] }
\eF
\[
+
\int_0^\tau \int_{\mathcal{T}^1}  \vr_R u^1_R \vu_R \cdot \vc{U} ( t, a, x_2 )   \ {\rm d}x_2 \ \dt
- \int_0^\tau \int_{\mathcal{T}^1}  \vr_L u^1_L \vu_L \cdot \vc{U} ( t, -a, x_2 )  \ {\rm d}x_2 \ \dt
\]
\[
+ \int_0^\tau \int_{\mathcal{T}^1} p(\vr_R) U^1 (t, a, x_2)\ {\rm d}x_2 \ \dt - \int_0^\tau \int_{\mathcal{T}^1}
p(\vr_L) U^1 (t, - a, x_2)\ {\rm d}x_2 \ \dt
\]
\[
= \int_0^\tau \intO{ \Big[ \vr \vu \cdot \partial_t \vc{U} + \vr [\vu \otimes \vu] : \Grad \vc{U} + p(\vr) \Div \vc{U} \Big] (t,x) } \ \dt.
\]

Next, we take $\varphi = \frac{1}{2} |\vc{U}|^2$  in (\ref{w1}):
\bFormula{w5}
\intO{ \left[ \frac{1}{2} \vr |\vc{U}|^2 (\tau, x) - \frac{1}{2} \vr_0 |\vc{U}|^2 (0, x) \right] }
\eF
\[
 +
\int_0^\tau \int_{\mathcal{T}^1}  \frac{1}{2} \vr_R u^1_R |\vc{U}|^2 ( t, a, x_2 )  \ {\rm d}x_2 \ \dt
- \int_0^\tau \int_{\mathcal{T}^1}  \frac{1}{2} \vr_L u^1_L |\vc{U}|^2 ( t , -a, x_2 )  \ {\rm d}x_2 \ \dt
\]
\[
= \int_0^\tau \intO{ \left[ \vr \vc{U} \cdot \partial_t \vc{U} + \vr \vu \cdot \Grad \vc{U} \cdot \vc{U} \right](t,x) } \ \dt.
\]

Summing up (\ref{w3}), (\ref{w4}), (\ref{w5}) we deduce
\bFormula{w6}
\intO{ \left[ \frac{1}{2} \vr |\vu -\vc{U}|^2 + H(\vr) \right](\tau, x) } - \intO{ \left[ \frac{1}{2} \vr_0 |\vu_0- \vc{U}(0, x)|^2 + H(\vr_0) \right] }
\eF
\[
+ \int_0^\tau \int_{\mathcal{T}^1} H(\vr_R) u^1_R \ {\rm d}x_2 \ \dt -
\int_0^\tau \int_{\mathcal{T}^1} H(\vr_L) u^1_L \ {\rm d}x_2 \ \dt
\]
\[
+ \int_0^\tau \int_{\mathcal{T}^1} \frac{1}{2} \vr_R |\vu_R - \vc{U}(t,a,x_2) |^2 u^1_R \ {\rm d}x_2 \ \dt - \int_0^\tau \int_{\mathcal{T}^1} \frac{1}{2} \vr_L |\vu_L - \vc{U}(t, -a, x_2) |^2 u^1_L \ {\rm d}x_2 \ \dt
\]
\[
+ \int_0^\tau \int_{\mathcal{T}^1} p(\vr_R) \left( u^1_R - U^1(t, a, x_2) \right) {\rm d}x_2 \ \dt - \int_0^\tau \int_{\mathcal{T}^1} p(\vr_L) \left( u^1_L - U^1(t, -a, x_2) \right) {\rm d}x_2 \ \dt
\]
\[
\leq
\int_0^\tau \intO{ \left[ \vr \left( \partial_t \vc{U} + \vu \cdot \Grad \vc{U} \right) \cdot (\vc{U} - \vu) - p(\vr) \Div \vc{U} \right](t,x) } \ \dt.
\]

\bigskip

{\bf Step 2}

Next step is to take $H'(r)$ as a test function in (\ref{w1}):
\bFormula{w7}
\intO{ \left[ \vr(\tau, x) H'(r) (\tau, x) - \vr_0 H'(r) (0, x) \right] }
\eF
\[
+
\int_0^\tau \int_{\mathcal{T}^1}  \vr_R u^1_R H'(r) ( t, a, x_2)  \ {\rm d}x_2 \ \dt
- \int_0^\tau \int_{\mathcal{T}^1}  \vr_L u^1_L H'(r) ( t, -a, x_2)  \ {\rm d}x_2 \ \dt
\]
\[
= \int_0^\tau \intO{ \left[ \vr \partial_t H'(r) + \vr \vu \cdot \Grad H'(r) \right] (t,x) } \ \dt;
\]
whence, in combination with (\ref{w6}), we deduce
\bFormula{w8}
\intO{ \left[ \frac{1}{2} \vr |\vu -\vc{U}|^2 + H(\vr) - H'(r) \vr \right](\tau, x) } - \intO{ \left[ \frac{1}{2} \vr_0 |\vu_0- \vc{U}(0, x)|^2 + H(\vr_0) - H'(r) (0, x) \vr_0 \right] }
\eF
\[
+ \int_0^\tau \int_{\mathcal{T}^1} \Big( H(\vr_R) - H'(r) (t, a, x_2) \vr_R \Big) u^1_R \ {\rm d}x_2 \ \dt -
\int_0^\tau \int_{\mathcal{T}^1} \Big( H(\vr_L)  - H'(r) (t, -a, x_2) \vr_L \Big)  u^1_L \ {\rm d}x_2 \ \dt
\]
\[
+ \int_0^\tau \int_{\mathcal{T}^1} \frac{1}{2} \vr_R |\vu_R - \vc{U}(t,a, x_2) |^2 u^1_R \ {\rm d}x_2 \ \dt - \int_0^\tau \int_0^1 \frac{1}{2} \vr_L |\vu_L - \vc{U}(t, -a, x_2) |^2 u^1_L \ {\rm d}x_2 \ \dt
\]
\[
+ \int_0^\tau \int_{\mathcal{T}^1} p(\vr_R) \left( u^1_R - U^1(t, a, x_2) \right) {\rm d}x_2 \ \dt - \int_0^\tau \int_{\mathcal{T}^1} p(\vr_L) \left( u^1_L - U^1(t, -a, x_2) \right) {\rm d}x_2 \ \dt
\]
\[
\leq \int_0^\tau \intO{ \left[ \vr \left( \partial_t \vc{U} + \vu \cdot \Grad \vc{U} \right) \cdot (\vc{U} - \vu) - p(\vr) \Div \vc{U} \right](t,x) } \ \dt
\]
\[
- \int_0^\tau \intO{ \left[ \vr \partial_t H'(r) + \vr \vu \cdot \Grad H'(r) \right](t,x) } \ \dt.
\]

\medskip

{\bf Step 3}

Finally, we write
\bFormula{w9}
\intO{ \left( H'(r) r - H(r) \right)(\tau, x) } - \intO{ \left( H'(r) r - H(r) \right) (0, x) } =
\int_0^\tau \intO{ \partial_t \left( H'(r) r - H(r) \right)(t,x) } \ \dt
\eF
\[
= \int_0^\tau \intO{ \partial_t p(r)(t,x) } \ \dt = \int_0^\tau \intO{ r \partial_t H'(r)(t,x) } \ \dt,
\]
and
\bFormula{w10}
0 = \int_0^\tau \intO{ \Div (p(r) \vc{U} )(t,x) } \ \dt - \int_0^\tau \int_{\mathcal{T}^1} p(r) U^1 (t, a, x_2) \ {\rm d}x_2 \ \dt+
\int_0^\tau \int_{\mathcal{T}^1} p(r) U^1 (t, -a, x_2) \ {\rm d}x_2
\eF
\[
= \int_0^\tau \intO{ \left( p(r) \Div \vc{U} + r \vc{U} \cdot \Grad H'(r) \right)(t,x)  } \ \dt
\]
\[
- \int_0^\tau \int_{\mathcal{T}^1} p(r) U^1 (t, a, x_2) \ {\rm d}x_2 \ \dt+
\int_0^\tau \int_{\mathcal{T}^1} p(r) U^1 (t, -a, x_2) \ {\rm d}x_2.
\]

Thus, combining (\ref{w8} - \ref{w10}), we obtain (\ref{REI}).

\qed

\section{Proof of Theorem \ref{Tm1}}
\label{p}

Having collected all the necessary material we are ready to prove Theorem \ref{Tm1}. The first observation is that the rarefaction wave solution
$[\tilde \vr, \tilde \vu]$ can be taken as test functions in the relative entropy inequality (\ref{REI}), specifically $r = \tilde \vr$, $\vc{U} =
\tilde \vu$. Indeed we have
\[
\vr, \tilde \vr, \vu, \tilde \vu \ \mbox{bounded}, \ \mbox{with} \ \partial_t \tilde \vr, \ \partial_t \tilde u^1 ,\
\partial_{x_1} \tilde \vr, \ \partial_{x_1} \tilde u^1 \in L^\infty(0,T; L^1 (\Omega));
\]
whence such a step may be justified via a density argument and the Lebesgue convergence theorem.

Accordingly, under the circumstances stated in Theorem \ref{Tm1}, the relative entropy inequality (\ref{REI}) simplifies to
\bFormula{REI1}
\intO{ \mathcal{E} \left( \vr, \vu \Big| \tilde \vr , \tilde \vu \right)(\tau,x) }
\eF
\[
 \leq \int_0^\tau \intO{ \left[ \vr \left( \partial_t \tilde u^1 + u^1  \partial_{x_1} \tilde u^1 \right)  (\tilde u^1 - u^1) + \Big( p(\tilde \vr) -  p(\vr) \Big) \partial_{x_1} \tilde u^1 \right](t,x) } \ \dt
\]
\[
+ \int_0^\tau \intO{ \Big[ (\tilde \vr - \vr) \partial_t H'(\tilde \vr) + (\tilde \vr \tilde u^1 - \vr u^1) \partial_{x_1} H'(\tilde \vr) \Big] (t,x) } \ \dt
\]
for a.a. $\tau \in (0,T)$.

Using the fact that $\tilde \vr$, $\tilde u^1$ satisfy the equations for $t > 0$ we deduce
\bFormula{s1}
\vr \left( \partial_t \tilde u^1 + u^1  \partial_{x_1} \tilde u^1 \right)  (\tilde u^1 - u^1) =
\vr \left( \partial_t \tilde u^1 + \tilde u^1  \partial_{x_1} \tilde u^1 \right)  (\tilde u^1 - u^1)
- \vr \partial_{x_1} \tilde u^1 (\tilde u^1 - u^1 )^2
\eF
\[
= - \frac{\vr}{\tilde \vr} \partial_{x_1} p(\tilde \vr) (\tilde u^1 - u^1) - \vr \partial_{x_1} \tilde u^1 (\tilde u^1 - u^1 )^2.
\]

Similarly,
\bFormula{s2}
\Big( p(\tilde \vr) -  p(\vr) \Big) \partial_{x_1} \tilde u^1 = -
p'(\tilde \vr) (\vr - \tilde \vr) \partial_{x_1} \tilde u^1 - \Big( p(\vr) - p'(\tilde \vr) (\vr - \tilde \vr) - p(\tilde \vr) \Big) \partial_{x_1} \tilde u^1,
\eF
and
\bFormula{s3}
(\tilde \vr - \vr) \partial_t H'(\tilde \vr) + (\tilde \vr \tilde u^1 - \vr u^1) \partial_{x_1} H'(\tilde \vr) =
p'(\tilde \vr) \partial_t \tilde \vr - \frac{\vr}{\tilde \vr} p'(\tilde \vr) \partial_t \tilde \vr + \tilde u^1 p'(\tilde \vr)
\partial_{x_1} \tilde \vr - \frac{\vr}{\tilde \vr} u^1 p'(\tilde \vr) \partial_{x_1} \tilde \vr.
\eF

Summing up (\ref{s1}-\ref{s3}) and using the fact that $\tilde \vr$, $\tilde u^1$ satisfy the equation of continuity (\ref{i1}) for $t > 0$, we may write
(\ref{REI1}) as
\bFormula{REI2}
\intO{ \mathcal{E} \left( \vr, \vu \Big| \tilde \vr , \tilde \vu \right)(\tau,x) }
\eF
\[
\leq - \int_0^\tau \intO{ \left[ \vr (\tilde u^1 - u^1)^2 + \Big( p(\vr) - p'(\tilde \vr) (\vr - \tilde \vr) - p(\tilde \vr) \Big)  \right]
\partial_{x_1} \tilde u^1 } \ \dt.
\]

Since the pressure $p$ is assumed to be convex, it is enough to observe that $\partial_{x_1} \tilde u^1 \geq 0$, in other words, the velocity component
of the rarefaction wave solution is non-decreasing in the spatial variable $x_1$. However, this follows easily from the standard analysis of the Riemann problem
that consists in rewriting the system in the Lagrangian coordinates and computing the solution in terms of the Riemann invariants. It turns out that there exist points $\xi^1_L \leq \xi^1_C \leq \xi^2_C \leq \xi^2_R$ such that
\bFormula{Rie}
u^1 (x^1/t) = \left\{ \begin{array}{c} u^1_L \ \mbox{for}\ x^1/t < \xi^1_L, \\ \\
R^1(x^1/t) \ \mbox{for}\ \xi^1_L \leq x^1/t \leq \xi^1_C, \\  \\
u^1_C \ \mbox{for} \ \xi^1_C < x^1/t < \xi^2_C, \\ \\
R^2(x^1/t) \ \mbox{for}\ \xi^2_C \leq x^1/t \leq \xi^2_R,\\ \\
u^1_R \ \mbox{for} \ x^1/t > \xi^2_R, \end{array} \right.
\eF
where $u^1_L \leq u^1_C \leq u^1_R$ are constants connected by the monotone (non-decreasing) functions $R^1$, $R^2$, see for instance \cite{ChanHsi}.

We have proved Theorem \ref{Tm1}.

\bRemark{conc}

The same method can be used to show uniqueness (in terms of the inital datum and the normal traces) of a general solution
$[\tilde \vr, \tilde \vu]$ in a space-time cylinder $(0,T) \times \Omega$, $\Omega \subset R^2$ as soon as:

\begin{itemize}
\item
$[\tilde \vr, \tilde \vu]$ is locally Lipschitz in the \emph{open} set $(0,T) \times \Omega$,
\item
$\tilde \vr \in W^{1,1}((0,T) \times \Omega)$, $\tilde \vu \in W^{1,1}((0,T) \times \Omega;R^2)$,
\item
$
\Grad \tilde\vu + \Grad^t \tilde \vu \geq - M \tn{I} \ \mbox{in}\ R^{2 \times 2}_{\rm sym} \ \mbox{a.a. in}\ (0,T) \times \Omega.
$
\end{itemize}
\eR

\bRemark{conc2}
Finally we summarize the current state of the art of the problem of (non)uniqueness of solutions to the Riemann problem for the compressible Euler system (\ref{i1}-\ref{i3}).
\begin{itemize}
\item If
\bFormula{Rdata1}
\left|\int_{\rho_L}^{\rho_R} \frac{\sqrt{p'(\tau)}}{\tau} {\rm d}\tau\right| \leq u_R^1 - u_L^1 < \int_0^{\rho_L} \frac{\sqrt{p'(\tau)}}{\tau} {\rm d}\tau + \int_0^{\rho_R} \frac{\sqrt{p'(\tau)}}{\tau} {\rm d}\tau
\eF
and $u_L^2 = u_R^2$, then the self-similar solution of the Riemann problem consists only of rarefaction waves and does not contain vacuum. We have proved in this paper that this solution is unique in the class of bounded entropy solutions.

\item If
\bFormula{Rdata2}
u_R^1 - u_L^1 < - \sqrt{\frac{(\rho_L-\rho_R)(p(\rho_L)-p(\rho_R))}{\rho_L\rho_R}}
\eF
and $u_L^2 = u_R^2$, then the self-similar solution of the Riemann problem consists of two shocks. Chiodaroli and Kreml proved in \cite{ChiKre} that for such data there exist infinitely many bounded entropy solutions.

\item Chiodaroli, De Lellis and Kreml proved in \cite{ChiDelKre} that there exist Riemann initial data such that $u_L^2 = u_R^2$ and
\bFormula{Rdata3}
 - \sqrt{\frac{(\rho_L-\rho_R)(p(\rho_L)-p(\rho_R))}{\rho_L\rho_R}} < u_R^1 - u_L^1 < \left|\int_{\rho_L}^{\rho_R} \frac{\sqrt{p'(\tau)}}{\tau} {\rm d}\tau\right|
\eF
(i.e. the self-similar solution consists of one shock and one rarefaction wave), for which there exists infinitely many bounded entropy solutions.

\end{itemize}
\eR

\def\cprime{$'$} \def\ocirc#1{\ifmmode\setbox0=\hbox{$#1$}\dimen0=\ht0
  \advance\dimen0 by1pt\rlap{\hbox to\wd0{\hss\raise\dimen0
  \hbox{\hskip.2em$\scriptscriptstyle\circ$}\hss}}#1\else {\accent"17 #1}\fi}


\end{document}